\documentclass[11pt]{article}
\usepackage{amsmath, graphicx, amsfonts,amssymb, calrsfs}
\usepackage{amsfonts,mathrsfs, color, amsthm}
\usepackage{bm}

\usepackage{enumitem}
\usepackage{hyperref}
\hypersetup{
  colorlinks   = true, %Colours links instead of ugly boxes
  urlcolor     = blue, %Colour for external hyperlinks
  linkcolor    = blue, %Colour of internal links
  citecolor   = red %Colour of citations
  }

\def\cK{\mathcal{K}}
\def\sphere{S^{n-1}}

\def\Rn{{\mathbb R^n}}
 \def\R{\mathbb{R}}

\newtheorem{theorem}{Theorem}[section]
\newtheorem{lemma}{Lemma}[section]
\newtheorem{remark}{Remark}[section]
\newtheorem{proposition}{Proposition}[section]
\newtheorem{corollary}{Corollary}[section]

\newtheorem{definition}{Definition}[section]

\def\bt{\begin{theorem}}
\def\et{\end{theorem}}
\def\bl{\begin{lemma}}
\def\el{\end{lemma}}
\def\br{\begin{remark}}
\def\er{\end{remark}}
\def\bc{\begin{corollary}}
\def\ec{\end{corollary}}
\def\bd{\begin{definition}}
\def\ed{\end{definition}}
\def\bp{\begin{proposition}}
\def\ep{\end{proposition}}
\def\ball{B^n_2}

\begin{document}
\title{The dual Orlicz-Minkowski problem }

\author{Baocheng Zhu, Sudan Xing and Deping Ye}
\date{}
\maketitle
\begin{abstract}

In this paper, the dual Orlicz curvature measure is proposed and its basic properties are provided. A variational formula for the dual Orlicz-quermassintegral is established in order to give a geometric interpretation of the dual Orlicz curvature measure.   Based on the established variational formula, a solution to the dual Orlicz-Minkowski problem regarding the dual Orlicz curvature measure is provided.

 \end{abstract}

 \section{Introduction}

   Given $p\in \R$ and a finite nonzero Borel measure $\mu$ defined on the unit sphere $\sphere\subset \Rn$, the $L_p$ Minkowski problem asks whether there exists  a convex body $K$ (i.e., a convex and compact subset in $\Rn$ with nonempty interior) such that $\mu$ is the $L_p$ surface area measure of $K$.   The $L_p$ Minkowski problem is arguably one of the most important problems in convex geometry. Like many other central objects of interest in convex geometry, such as the $L_p$ affine and geominimal surface areas and the $L_p$ John ellipsoids (see e.g., \cite{John1948, Lu96, LYZ05,  Petty1974, Petty1985, Ye2015a, Zhao-IMRN}), the $L_p$ Minkowski problem for $p\geq 1$ is closely related to optimization problems involving the $L_p$ mixed volume of convex bodies; this can be seen intuitively from the equality characterization of the $L_p$ Minkowski inequality for the $L_p$ mixed volume. When $p=1$, the $L_p$ Minkowski problem becomes the classical Minkowski problem which went back to Minkowski \cite{min1897,min1903}. The $L_p$ Minkowski problem was first posed by Lutwak \cite{Lu93}. Since then, the $L_p$ Minkowski problem attracted a lot of attention and amazing progress has been made (see e.g., \cite{chen, chouw06,HuMaShen, HugLYZ,JLZ2016,LYZ04,Umanskiy, zhug2015a,zhug2015b,zhug2017}). When $p=0$,  it becomes  the logarithmic Minkowski problem (see e.g.,
\cite{BHZ2016,BLYZ2013,stancu02,stancu03,stancu08,zhug2014}).  Note that solutions to the $L_p$ Minkowski problem are the key ingredients in the rapidly developing $L_p$ Brunn-Minkowski theory of convex bodies. For instance, it has been used to establish the sharp affine Sobolev inequalities \cite{CLYZ-1, HS-JFA-1, LYZ02, zh99b}.

 Recently, much effort has been made to develop the nonhomogeneous theory analogous to the $L_p$ Brunn-Minkowski theory; such a new theory is called the Orlicz-Brunn-Minkowski theory of convex bodies. This new theory started from the works of  Lutwak, Yang and Zhang \cite {LYZ2010a,LYZ2010b} and Ludwig \cite{Ludwig2010}, and was greatly pushed forward by Gardner, Hug and Weil \cite{GHW2014} and independently by Xi, Jin and Leng \cite{XJL}, due to the discovery of the Orlicz addition of convex bodies. In particular, the Orlicz mixed volumes of convex bodies have been obtained in \cite{GHW2014, XJL}, which were used to discover, for instance, the Orlicz affine and geominimal surface areas and the Orlicz-John ellipsoid  \cite{Ye2015a2, YJL-2015, ZHY2017, zouxiong}. On the other hand, the Orlicz-Minkowski problem can be asked with the $L_p$ surface area measure (in the $L_p$ Minkowski problem) replaced by the Orlicz surface area measure. Solutions to the Orlicz-Minkowski problem can be found in \cite{HLYZ2010,huanghe2012,liaijun2014}.

Replacing convex bodies and their $L_p$ addition by star bodies and their ($L_p$) radial addition, Lutwak
\cite{Lu75, Lu88} developed the beautiful dual ($L_p$) Brunn-Minkowski theory for star bodies. An amazing success, among others,  of the dual ($L_p$) Brunn-Minkowski theory is to provide powerful tools to solve the famous Busemann-Petty problem (see e.g., \cite{G1,G2,G3, GardnerKoldobski1999, Koldobski1997, zh99a}).  Many notions in the $L_p$ Brunn-Minkowski theory for convex bodies have their dual analogues in the dual theory (see e.g., \cite{Lu75, Lu88} or the book \cite{Gard} written by Gardner for more background and references). However, it is only very recent that the dual curvature measures, which are dual to the surface area measures, were discovered in the seminal work \cite{HLYZ} by Huang, Lutwak, Yang and Zhang.  Regarding the dual curvature measures, they posed the following dual Minkowski problem for the $q$-th dual curvature measure: \emph{for $q\in \R$ and $\mu$ a given nonzero finite Borel measure on $S^{n-1}$, can we find a convex body (ideally with the origin in its interior) such that $\mu$ is equal to the $q$-th dual curvature measure of $K$?} When $q=n$, the dual Minkowski problem becomes the logarithmic Minkowski problem. The authors in  \cite{HLYZ} provided the existence of solutions (i.e., origin-symmetric convex bodies) to the dual Minkowski problem with $q\in (0,n)$ for even measure $\mu$. For $q<0$, the existence and uniqueness of the solution to the dual Minkowski problem were given recently in \cite{zhao} by Zhao. See also \cite{BHP, zhao-jdg} for more works on the dual curvature measure and the related dual Minkowski problem.

It is our main goal in this paper to develop the dual Orlicz curvature measure and provide solutions to the dual Orlicz-Minkowski problem. The dual Orlicz curvature measure will be defined similarly to the $q$-th dual curvature measure, which relies on the linear Orlicz radial addition of convex bodies (indeed of star bodies) discovered by Gardner, Hug, Weil and Ye \cite{ghwy15} and independently by Zhu, Zhou and Xu \cite{Zhub2014}. In fact, with the help of the linear Orlicz radial addition, the authors in \cite{ghwy15, Zhub2014} obtained the formulas for the dual Orlicz mixed volume of convex bodies (indeed of star bodies). Note that the dual Orlicz mixed volume plays fundamental role in the dual Orlicz-Brunn-Minkowski theory for star bodies developed in \cite{ghwy15, Zhub2014}. In particular, it is crucial for developing the theory of the dual Orlicz affine and geominimal surface areas \cite{Ye2016a} and the  Orlicz-Legendre ellipsoids \cite{zouxiong2016}.  See e.g., \cite{JYL, zhaoc, zhug2012} for more works in the rapidly developing dual Orlicz-Brunn-Minkowski theory.

 Let $K$ be  a convex body  with the origin in its interior and   $\varphi: (0, \infty)\rightarrow (0, \infty)$ be a continuous function. We propose the following definition for the dual Orlicz curvature measure $\widetilde{C}_\varphi(K,\cdot)$: for each Borel set $\eta\subset \sphere$,  let
 $$
 \widetilde{C}_\varphi(K,\eta)
 =\frac{1}{n}\int_{\pmb{\alpha}^*_K(\eta)}\varphi(\rho_K(u))du,
 $$ where $\rho_K$ is the radial function of $K$,  $\pmb{\alpha}^*_K$ is the reverse radial Gauss image on $S^{n-1}$,  and $\,du$ is the spherical measure on $\sphere$  (see Section 2 for definitions and notations). The dual Orlicz $L_{\varphi}$ quermassintegral of $K$, denoted by $\widetilde{V}_\varphi(K),$  is just   $$\widetilde{V}_\varphi(K)=
 \widetilde{C}_\varphi(K,\sphere)
 =\frac{1}{n}\int_{\sphere}\varphi(\rho_K(u))du.
 $$ When $\varphi(t)=t^q$ with $q\in\R$, one gets the $q$-th dual curvature measure and the $q$-th dual quermassintegral of $K$  \cite{HLYZ}. We are interested in the following dual Orlicz-Minkowski problem: {\it under what conditions on $\varphi$ and a given nonzero finite Borel measure $\mu$ on $\sphere$,  there exist a constant $\tau>0$ and a convex body $K$ (ideally with the origin in its interior) such that $\mu=\tau \widetilde{C}_\varphi(K,\cdot)$?} When $\varphi(t)=t^q$ for $0\neq q\in \R$, this becomes the  dual Minkowski problem for the $q$-th dual curvature measure posed in \cite{HLYZ}.

  Suppose that $\varphi$ (and its companion $\phi$) satisfies conditions A1)-A3) (see the detailed description of these conditions in Section \ref{variation-2--2--1-1}). Our solution to the dual Orlicz-Minkowski problem is stated  and proved in Theorem \ref{solution-dual-Orlicz-main}. Note that the condition on $\mu$ in Theorem \ref{solution-dual-Orlicz-main} is the minimal requirement for the solutions to various Minkowski problems. For a nonzero finite Borel measure $\mu$ on $\sphere$, let  $$|\mu|=\int_{\sphere}\,d\mu.$$ Let $a_+=\max\{a, 0\}$ for each $a\in \R$.

 \vskip 2mm \noindent  {\bf Theorem \ref{solution-dual-Orlicz-main}.} {\it Let $\mu$ be a nonzero finite Borel measure on $S^{n-1}$. Then the following  statements are equivalent:
\begin{itemize} \item [i)] $\mu$ is not concentrated in any closed hemisphere, i.e.,   $$\int _{\sphere} \langle \xi, \theta\rangle_+ \,d\mu(\theta)>0\ \ \ \mathrm{for\ \  all\ \ } \xi\in \sphere;$$
\item [ii)]  there exists a convex body $K$ with the origin in its interior, such that, $$\frac{\mu}{|\mu|}=\frac{\widetilde{C}_\varphi(K,\cdot)}{\widetilde{V}_\varphi(K)}. $$   \end{itemize}}

Our paper is organized as follows. Section \ref{section 2} dedicates to the necessary notations and background in order to present our results. The definition and properties for the dual Orlicz curvature measure are given in Section \ref{dual Orlicz curvature measures---1--1}. Section  \ref{variation-2--2--1-1} aims to establish a variational formula for the dual Orlicz-quermassintegral, which provides a geometric interpretation for the dual Orlicz curvature measure. Such a variational formula is the key for our solution to the
dual Orlicz-Minkowski problem.  Finally, in Section \ref{solution-dual-Orlicz}, we present a solution to the dual Orlicz-Minkowski problem, i.e., Theorem  \ref{solution-dual-Orlicz-main}.

\section{Background and Notations}\label{section 2}
Basic terminologies and well-known facts in convex geometry required for presenting our results are collected. When talking about the concepts for the dual curvature measures, we manage to keep our  notations
as consistent as possible with  those in \cite{HLYZ}.

In $\mathbb{R}^{n}$, the standard inner product and the usual Euclidean norm are denoted by  $x\cdot y$ and $|x|$, respectively, for $x, y\in \mathbb{R}^{n}$. The set $\ball$ refers to the unit Euclidean ball in $\Rn$ and its boundary, i.e., the unit sphere in $\Rn$, is denoted by $\sphere$. The origin in $\Rn$ is denoted by $o$. For each $x\in\R^n\setminus \{o\}$, let $\bar{x}=x/|x|\in \sphere$ be the direction vector of $x$. For  a set $E\subset \Rn$, $\partial E$ means the boundary of $E$ and $V(E)$ denotes the volume of $E$ (if it exists).

A set $E\subset \Rn$ is convex if $\lambda x+(1-\lambda)y\in E$ for all $x, y\in E$ and $\lambda \in [0, 1]$.  For a compact and convex set $K\subset\R^n$, one can define a sublinear function with positive homogeneity of degree $1$, $h_K: S^{n-1}\rightarrow \R,$ by $$
h_K(u)=\max\big\{x\cdot u: x\in K\big\}, \ \ \text{for each}\ \ u\in S^{n-1}.
$$ The function $h_K$ is called the support function of $K$.  Define the Hausdorff distance between two convex and compact sets $K, L\subset \Rn$ by $$d_H(K, L)=\|h_K-h_L\|_{\infty}=\sup_{u\in \sphere} |h_K(u)-h_L(u)|.$$ We say a sequence of convex and compact sets $\{K_i\}_{i=1}^{\infty}$ in $\Rn$ converges to a convex and compact set $K\subset \Rn$ with respect to the Hausdorff metric, if $\lim_{i\rightarrow \infty} d_H(K_i, K)=0.$ A convex body is a convex and compact set $K\subset \Rn$ with nonempty interior. In particular, we work on the set $\mathcal{K}_0^n$ which contains all convex bodies with the origin $o$ in their interiors and hence $h_K>0$ for each $K\in \mathcal{K}_0^n$.

Associated to each convex body $K\in \mathcal{K}_0^n$, one can also define the so-called  radial
function $\rho_K: S^{n-1}\rightarrow \R$  by
$$
\rho_K(u)=\max\big\{\lambda:\ \lambda u\in K\big\}, \ \ \text{for each}\ \ u\in S^{n-1}.
$$ It is easily checked that $\rho_K(u)u\in \partial K$ for all $u\in \sphere$.  Similar to the Hausdorff metric, one can define the radial  metric $d_{\rho}(\cdot, \cdot)$ on $\cK_0^n$; namely for two convex bodies $K, L\in \cK_0^n$,
$$d_{\rho}(K, L)=\|\rho_K-\rho_L\|_{\infty}=\sup_{u\in \sphere} |\rho_K(u)-\rho_L(u)|.$$ Note that, on $\cK_0^n$, the Hausdorff metric is equivalent to the radial metric. Two convex bodies $K, L\in \cK_0^n$ are dilates of each other if and only if $h_K=\lambda h_L$ (or $\rho_K=\lambda \rho_L$) on $\sphere$ for some constant $\lambda>0$. 

For each  $K\in \cK_0^n$, the set $$
K^*=\{x\in \R^n: x\cdot y\le 1 \ \ \text{for all}\ \ y\in K\}
$$ defines a convex body, called the polar body of $K$, and clearly $K^*\in  \cK_0^n$. Regarding the support and radial functions of $K$ and $K^*$, one has  $$
\rho_K\cdot h_{K^*}=1\ \ \ \mathrm{and}\ \ \ h_K\cdot \rho_{K^*}=1\ \ \ \mathrm{on}\ \ \sphere.$$
The bipolar theorem asserts that $(K^*)^*=K$ for each $K\in \cK_0^n$. Moreover, if $K_i, K\in \cK_0^n$ such that $K_i\rightarrow K$ with respect to the Hausdorff metric, then $K_i^* \rightarrow K^*$, as well, with respect to the Hausdorff metric.

Consider the space of continuous functions on $S^{n-1}$, denoted by $C(S^{n-1})$, with the maximal norm $\|\cdot\|_{\infty}$. The convergence $f_i\rightarrow f$ with $f_i, f \in C(S^{n-1})$ is in the sense of $$\lim_{i\rightarrow \infty}  \sup_{u\in \sphere} |f_i(u)-f(u)|=\lim_{i\rightarrow \infty} \|f_i-f\|_{\infty}=0.$$  The subset of positive functions in $C(S^{n-1})$  is denoted by $C^+(S^{n-1})$.   The set $\Omega\subset S^{n-1}$ is always assumed to be a closed set that is not contained in any closed hemisphere of $S^{n-1}$. The Wulff shape associated to a continuous function $f: \Omega\rightarrow (0,\infty)$, denoted by $[f]$, is the convex body with $o$ in its interior such that
$$
[f]=\cap_{u\in \Omega}\{x\in \Rn: x\cdot u\leq f(u)\}.
$$ Clearly $h_{[f]}\le f$ if $\Omega=S^{n-1}$. Moreover, when $f=h_K$ is the support function of a convex body $K\in \cK_0^n$, one has $[h_K]=K$. The convex hull $\langle \rho \rangle$ generated by a continuous function $\rho: \Omega\rightarrow (0,\infty)$ is formulated by
$$
\langle \rho \rangle=\mathrm{conv}\big\{\rho(u)u: u\in \Omega\big\}.
$$ It can be checked that $\langle \rho \rangle\in \cK_0^n$ (hence  $h_{\langle \rho \rangle}>0$), again since $\Omega$ is not contained in any closed hemisphere of $S^{n-1}$.  Moreover, when $\rho=\rho_K$ is the radial function of $K\in \cK_0^n$, then $\langle \rho_K \rangle=K$.

The following lemma, i.e., \cite[Lemma 2.8]{HLYZ}, is about the relation between the Wulff shape and the convex hull.

 \bl \label{relation}
 Let $\Omega\subset S^{n-1}$ be a closed set that is not contained in any closed hemisphere of $S^{n-1}$. For each continuous function $f: \Omega \rightarrow (0,\infty)$, one has
 $$
 [f]^*=\langle 1/f \rangle.
 $$
 \el

 The convergence of the Wulff shapes can be obtained by the convergence of continuous functions, according to the Aleksandrov's convergence theorem \cite{Aleks1938,Sch}: if $f_i, f: \Omega \rightarrow (0,\infty)$ are continuous functions for all $i\geq 1$ such that $\max_{u\in \Omega}|f_i(u)-f(u)|\rightarrow 0,$ then $[f_i]\rightarrow [f]$ with respect to the Hausdorff metric. Together with Lemma \ref{relation}, the convergence of convex hulls can also be obtained by the convergence of functions: if $\rho_i, \rho: \Omega \rightarrow (0,\infty)$  are continuous functions for all $i\geq 1$  such that $\max_{u\in \Omega} |\rho_i(u)-\rho(u)|\rightarrow 0,$ then $\langle \rho_i \rangle \rightarrow \langle \rho \rangle$ with respect to the Hausdorff metric.

Let $K\in\cK_0^n$ and $u\in S^{n-1}$. The set $
H(K,u)=\big\{x\in\R^n: x\cdot u=h_K(u)\big\}$ defines the supporting hyperplane of $K$ at the direction $u$. For each $\sigma \subset \partial K$, let $\pmb{\nu}_K(\sigma)\subset \sphere$ be the spherical image of $\sigma$, i.e.,
$$
\pmb{\nu}_K(\sigma)=\big\{u\in S^{n-1}: x\in H(K,u)\ \ \text{for some}\ \ x\in \sigma\big\}.
$$ The reverse spherical image  $\pmb{x}_K(\eta)$  of each $\eta\subset S^{n-1}$ is a subset of $\partial K$ defined by  $$
\pmb{x}_K(\eta)=\big\{x\in\partial K: x\in H(K,u)\ \ \text{for some}\ \ u\in \eta\big\}.
$$ Let $\omega$ be a subset of $S^{n-1}$. The subset  $\pmb{\alpha}_K(\omega)\subset \sphere$ given by
$$
\pmb{\alpha}_K(\omega)=\pmb{\nu}_K(\{\rho_K(u)u\in\partial K: u\in \omega\})
$$ is called the radial Gauss image  of $\omega$.   We often write $\pmb{\nu}_K(\{x\})$, $\pmb{x}_K(\{u\})$ and  $\pmb{\alpha}_K(\{u\})$ by $\pmb{\nu}_K(x)$, $\pmb{x}_K(u)$ and  $\pmb{\alpha}_K(u)$, respectively. When $\pmb{\nu}_K(x)$, $\pmb{x}_K(u)$ and $\pmb{\alpha}_K(u)$ all contain only one element, they will be written as  $\nu_K(x)$, $x_K(u)$ and $\alpha_K(u)$, respectively.

 For $K\in \cK_0^n$, define $\sigma_K\subset \partial K$, $\eta_K\subset \sphere$ and  $\omega_K\subset S^{n-1}$, respectively,  by
 \begin{eqnarray*} \sigma_K&=&\big\{x\in\partial K: \ \ \  \pmb{\nu}_K(x) \ \ \mathrm{has\ more\ than\ one\ element}\big\};\\
 \eta_K&=& \big\{v \in \sphere: \  \pmb{x}_K(v) \ \ \mathrm{has\ more\ than\ one\ element}\big\};\\ \omega_K&=& \big\{v \in \sphere: \  \pmb{\alpha}_K(v) \ \ \mathrm{has\ more\ than\ one\ element}\big\}. \end{eqnarray*} It is well-known that  $\mathcal{H}^{n-1}(\sigma_K)=0$ according to \cite[p.84]{Sch}, and $\eta_K$ and $\omega_K$  have  spherical measure zero according to \cite[Theorems 2.2.5 and 2.2.11]{Sch} or \cite[p.339-340]{HLYZ}. Hereafter,  the standard notation $\mathcal{H}^{n-1}_K$ of $K\in \cK_0^n$,  more often abbreviated by $\mathcal{H}^{n-1}$, is for the $(n-1)$-dimensional Hausdorff measure on $\partial K$. For all $u\in S^{n-1}\setminus\omega_K$, one sees that
$\alpha_K(u)=\nu_K\circ(\rho_K(u)u).$  In particular, $
\alpha_K(\bar{x})=\nu_K(x)
$ for all $x\in \partial K$ such that $\bar{x}\in S^{n-1}\setminus\omega_K$.  Note that $\nu_K(x)$, $x_K(u)$ and $\alpha_K(u)$ are all continuous \cite{HLYZ, Sch}.

For $K\in \cK_0^n$, let $\partial^\prime K=\partial K\setminus\sigma_K$. This implies that $\mathcal{H}^{n-1}(\partial^\prime K)>0$. Associated to each $K\in \cK_0^n$, the surface area measure $S(K,\cdot)$  defined on $\sphere$ is the measure with the following property: for any Borel set $\eta\subset \sphere$,
$$
S(K,\eta)=\mathcal{H}^{n-1}(\nu_K^{-1}(\eta)).
$$ Moreover, for each continuous function $g:S^{n-1}\rightarrow \R$, one has (see e.g., \cite[(2.12)]{HLYZ})
 \begin{equation}\label{surface area measure}
 \int_{S^{n-1}}g(u)dS(K,u)=\int_{\partial^\prime K} g(\nu_K(x))d\mathcal{H}^{n-1}(x).
 \end{equation}

  For  each subset $\eta\subset S^{n-1}$, let $$
\pmb{\alpha}^*_K(\eta)=\big\{\bar{x}: x\in \partial K \cap H(K,u)\ \text{for some}\ u\in \eta\big\}.
$$ The set $\pmb{\alpha}^*_K(\eta)$ is called the reverse radial Gauss image of $\eta.$ Note that
$$
\pmb{\alpha}^*_K(\eta)=\overline{\pmb{x}_K(\eta)}\subset S^{n-1}.
$$
 Let $\alpha^*_K(u)= \overline{x_K(u)}$ for each $u\in S^{n-1}\setminus\eta_K$, and  $\alpha^*_K$ is called the reverse radial Gauss map of $K$. Note that
 $\alpha^*_K$ is continuous. For each $\eta \subset S^{n-1}$ and for almost all $u\in S^{n-1}$ with respect to the spherical measure, one has  (see \cite[(2.21)]{HLYZ})    \begin{equation}\label{map-reverse}
u\in \pmb{\alpha}_K^*(\eta)\ \ \ \text{if and only if}\ \ \ \alpha_K(u)\in \eta.
 \end{equation}

The following lemma is the combination of \cite[Lemmas 2.5 and 2.6]{HLYZ}.
 \bl \label{relation-radial Gauss image}
 Let $K\in \cK_0^n$. For each $\eta\subset S^{n-1}$, one has $
 \pmb{\alpha}^*_K(\eta)=\pmb{\alpha}_{K^*}(\eta).$
 Moreover, for almost all $v\in S^{n-1}$ with respect to the spherical measure, one has $\alpha^*_K(v)=\alpha_{K^*}(v).$
 \el

The following lemma is the combination of \cite[Lemmas 2.1-2.4]{HLYZ}.
\bl   \label{properties of maps} Let $K\in \cK_0^n$ be a convex body with $o$ in its interior. \begin{itemize} \item [i)]  If $\eta\subset S^{n-1}$ is a Borel set, then $\pmb{\alpha}_K^*(\eta)=\overline{\pmb{x}_K(\eta)}\subset S^{n-1}$ is spherical measurable.
\item [ii)]  Let $K_i\in\cK_0^n$ be such that $\lim_{i\rightarrow \infty}K_i=K_0\in\cK_0^n$ with respect to the Hausdorff metric. Let $\omega=\cup_{i=0}^\infty \omega_{K_i}$ be the set (of spherical measure zero) of which all of the $\alpha_{K_i}$ are defined. If $u_i\in S^{n-1}\setminus \omega$ are such that $\lim_{i\rightarrow \infty}u_i=u_0\in S^{n-1}\setminus \omega$, then $\lim_{i\rightarrow \infty}\alpha_{K_i}(u_i)=\alpha_{K_0}(u_0).$
\item [iii)]  If $\{\eta_j\}_{j=1}^\infty$ is a sequence of subsets of $S^{n-1}$, then
$
\pmb{\alpha}^*_K\big(\!\cup_{j=1}^\infty\eta_j\big)
=\cup_{j=1}^\infty\pmb{\alpha}^*_K(\eta_j).
$
\item [iv)]  If $\{\eta_j\}_{j=1}^\infty$ is a sequence of pairwise disjoint sets in $S^{n-1}$, then $\{\pmb{\alpha}_K^*(\eta_j)\setminus \omega_K\}_{j=1}^\infty$ is pairwise disjoint as well.\end{itemize}
\el
 
 For more background in convex geometry, in particular the notions related to the radial Gauss map, please see   \cite{Gard, Gruber2007, HLYZ, Sch}.

\section{The dual Orlicz curvature measure}\label{dual Orlicz curvature
measures---1--1}

Let $\varphi_i: (0,
\infty)\rightarrow (0, \infty)$ be  strictly increasing continuous functions with $\lim_{t\rightarrow 0^+} \varphi_i(t)=0$ and $\lim_{t\rightarrow \infty} \varphi_i(t)=\infty$, $i=1, 2$.
For
$\varepsilon>0$ and two convex bodies $K, L\in \mathcal{K}_0^n$, define $\rho_{K\widetilde{+}_{\varphi, \varepsilon}L}: \sphere\rightarrow \R$ by \cite{ghwy15, Zhub2014} \begin{equation}
\varphi_1\bigg(\frac{\rho_K(u)}{\rho_{K\widetilde{+}_{\varphi, \varepsilon}L}(u)}\bigg)+
\varepsilon\varphi_2\bigg(\frac{\rho_L(u)}{\rho_{K\widetilde{+}_{\varphi, \varepsilon}L}(u)}\bigg)=1\
\ \ \mbox{for} \ \ u\in S^{n-1}.\label{dual-orlicz-addition-11}
\end{equation} Clearly $\rho_{K\widetilde{+}_{\varphi, \varepsilon}L}$ is a continuous function on $\sphere$ and  $K\widetilde{+}_{\varphi,
\varepsilon}L$ is called the linear Orlicz radial addition of $K, L\in \mathcal{K}_0^n$. When $\varphi_i: (0,
\infty)\rightarrow (0, \infty)$  are strictly decreasing continuous with $\lim_{t\rightarrow 0^+} \varphi_i(t)=\infty$ and $\lim_{t\rightarrow \infty} \varphi_i(t)=0$, $i=1, 2$, the function $\rho_{K\widetilde{+}_{\varphi, \varepsilon}L}$ can also be defined by formula (\ref{dual-orlicz-addition-11}). If $\varphi_1'(1)$, the  derivative of $\varphi_1$ at $1$, exists and is nonzero,  the following variational formula holds \cite{ghwy15, Zhub2014} :
$$ \varphi_1'(1)\lim_{\varepsilon\rightarrow 0^+}\frac{V(K\widetilde{+}_{\varphi, \varepsilon}L)-V(K)}{\varepsilon}
=\int_{S^{n-1}}\varphi_2\left(\frac{\rho_L(u)}{\rho_K(u)}\right)\rho_K(u)^n\,du.
$$ Motivated by this formula, one can define the dual Orlicz mixed volume of $K, L\in \cK_0^n$ by
 \begin{equation*}
\widetilde{V}_{\psi}(K, L)
=\frac{1}{n}\int_{S^{n-1}}\psi\left(\frac{\rho_L(u)}{\rho_K(u)}\right)\rho_K(u)^n\,du,
 \end{equation*} where $\psi: (0, \infty)\rightarrow (0, \infty)$ is a  continuous function.  In particular,  $$\widetilde{V}_\psi(K, \ball)
=\frac{1}{n}\int_{S^{n-1}}\psi\left(\frac{1}{\rho_K(u)}\right)\rho_K(u)^n\,du.
$$ The above definitions and results are stated here only for convex bodies, however they also hold for more general star sets, see details in \cite{ghwy15, Zhub2014}.

Let $\phi: (0, \infty)\rightarrow (0, \infty)$ be a  continuous function such that  $$\phi(\rho_K)=\psi\left(\frac{1}{\rho_K}\right)\rho_K^n.$$ We now propose the following definition for the dual Orlicz-quermassintegral of  $K\in \cK_0^n$.

 \bd\label{definition-volume}
Let $K\in\cK_0^n$ and $\phi: (0,\infty)\rightarrow (0,\infty)$ be a continuous function. Define the dual Orlicz-quermassintegral $\widetilde{V}_\phi(K)$  by
 \begin{equation}\label{formula-volume}
 \widetilde{V}_\phi(K)=\frac{1}{n}\int_{S^{n-1}}\phi(\rho_K(u))du.
 \end{equation}
 \ed

The continuity of $\widetilde{V}_\phi$ on $\cK_0^n$ is stated in the following lemma.

 \bl\label{continuity-dual-qu-1}
 Let $\phi: (0,\infty)\rightarrow (0,\infty)$ be a continuous function. Suppose that the sequence $\{K_i\}_{i=1}^\infty\subset \cK_0^n$ converges to $K\in\cK_0^n$ with respect to the Hausdorff metric. Then
 $$
 \lim_{i\rightarrow \infty}\widetilde{V}_\phi(K_i)=\widetilde{V}_\phi(K).
 $$
 \el
 \begin{proof} The convergence $K_i\rightarrow K$ with respect to the Hausdorff metric  implies that $\rho_{K_i}(u)\rightarrow \rho_{K}(u)$ uniformly on $S^{n-1}$. As $K\in \cK_0^n$,  there exist constants $r, R>0$, such that, for all $u\in \sphere$ and for all $i=1, 2, \cdots$,
$$
r\ \le \ \rho_{K_i}(u), \rho_{K}(u) \ \le R.
$$ It follows from the continuity of $\phi$ on $[r, R]$ that
$\phi(\rho_{K_i})\le C$ for some constant $C>0$. The desired continuity follows immediately from  the dominated convergence theorem:
$$\lim_{i\rightarrow \infty}\widetilde{V}_\phi(K_i)=\lim_{i\rightarrow \infty} \frac{1}{n} \int_{S^{n-1}} \phi(\rho_{K_i}(u))\,du=\frac{1}{n}\int_{S^{n-1}} \phi(\rho_{K}(u))\,du=\widetilde{V}_\phi(K).
$$ \end{proof} The dual Orlicz-quermassintegral $\widetilde{V}_\phi(\cdot)$ for convex bodies can be adopted to define its analogue for functions through the Wulff shape. That is, for each $h\in C^+(S^{n-1})$, we define $\widetilde{V}_\phi([h])$ to be the dual Orlicz-quermassintegral of $h$. Lemma \ref{continuity-dual-qu-1} and the Aleksandrov's convergence theorem yield the continuity of  the dual Orlicz-quermassintegral  on $C^+(S^{n-1})$. That is, if $h_i\rightarrow h$ with $h_i, h\in C^+(S^{n-1})$, then
  $$
  \lim_{i\rightarrow \infty}\widetilde{V}_\phi([h_i])=\widetilde{V}_\phi([h]).
  $$

The dual Orlicz  curvature measure is defined as follows.
 \bd\label{definition-general curvature measure}
Let $K\in\cK_0^n$ and $\varphi: (0,\infty)\rightarrow (0,\infty)$ be a continuous function. The dual Orlicz curvature measure of $K$, denoted by $\widetilde{C}_\varphi(K,\cdot)$, is defined to be the measure such that for each Borel set $\eta\subset S^{n-1}$,
 \begin{equation}\label{formula-general curvature measure}
 \widetilde{C}_\varphi(K,\eta)=\frac{1}{n}\int_{\pmb{\alpha}^*_K(\eta)}\varphi(\rho_K(u))\,du
 =\frac{1}{n}\int_{S^{n-1}}\emph{\textbf{1}}_{\pmb{\alpha}^*_K(\eta)}(u)\varphi(\rho_K(u))\,du.
 \end{equation} \ed Clearly, $ \widetilde{V}_\varphi(K)=\widetilde{C}_\varphi(K,\sphere)$. When $\varphi(t)=t^q$ with $q\in \R$, for each $K\in\cK_0^n$, one gets the $q$-th dual curvature measure  $\widetilde{C}_q(K,\cdot)$  \cite[Definition 3.2]{HLYZ}: for each Borel set $\eta\subset \sphere$,
$$
 \widetilde{C}_q(K,\eta)=\frac{1}{n}\int_{\pmb{\alpha}^*_K(\eta)} \rho_K^q(u) \,du
 =\frac{1}{n}\int_{S^{n-1}}\mathrm{\textbf{1}}_{\pmb{\alpha}^*_K(\eta)}(u)  \rho_K^q(u)\,du.
$$  In Proposition \ref{measure-check}, we will prove that $\widetilde{C}_\varphi(K,\cdot)$ is indeed a  Borel measure on $S^{n-1}$ after we prove the following useful result. When $\varphi(t)=t^q$ for $q\in \R$, one gets \cite[Lemma 3.3]{HLYZ}.

  \bl \label{measure-change}
 Let $K\in\cK_0^n$ and $\varphi: (0,\infty)\rightarrow (0,\infty)$ be a continuous function. For each bounded Borel function $g: S^{n-1}\rightarrow \R$, one has
 \begin{equation}\label{integral-general curvature measure}
 \int_{S^{n-1}}g(v)d\widetilde{C}_\varphi(K,v)=
 \frac{1}{n}\int_{S^{n-1}}g(\alpha_K(u))\varphi(\rho_K(u))du.
 \end{equation}
 \el

\begin{proof} Let $g: S^{n-1}\rightarrow \R$ be a bounded Borel function. As explained in \cite[p.353]{HLYZ}, both $g$ and $g\circ \alpha_K$ are Lebesgue integrable on $\sphere$. Since $g$ is bounded, the desired formula (\ref{integral-general curvature measure}) follows, by the dominated convergence theorem, if formula (\ref{integral-general curvature measure}) is proved for simple functions. Consider the simple function $
 \gamma=\sum_{i=1}^m c_i \textbf{1}_{\eta_i}$  with $c_i\in\R$ and Borel sets  $\eta_i\subset S^{n-1}$. Then \begin{eqnarray*}
 \int_{S^{n-1}}\gamma(v)d\widetilde{C}_\varphi(K,v)
 =\int_{S^{n-1}}\sum_{i=1}^m c_i \textbf{1}_{\eta_i}(v)d\widetilde{C}_\varphi(K,v)
 =\sum_{i=1}^m c_i \widetilde{C}_\varphi(K,\eta_i).
 \end{eqnarray*} It follows from Definition \ref{definition-general curvature measure} and (\ref{map-reverse}) that
 \begin{eqnarray*}
 \int_{S^{n-1}}\gamma(v)d\widetilde{C}_\varphi(K,v)
 &=&\frac{1}{n}\int_{S^{n-1}}\sum_{i=1}^m c_i \textbf{1}_{\pmb{\alpha}^*_K(\eta_i)}(u)\varphi(\rho_K(u))du\\
 &=&\frac{1}{n}\int_{S^{n-1}}\sum_{i=1}^m c_i \textbf{1}_{\eta_i}(\alpha_K(u))\varphi(\rho_K(u))du\\
 &=&\frac{1}{n}\int_{S^{n-1}} \gamma(\alpha_K(u)) \varphi(\rho_K(u))du.
 \end{eqnarray*}
 This proves (\ref{integral-general curvature measure}) for simple functions and hence for all bounded Borel functions $g: S^{n-1}\rightarrow \R$.
 \end{proof}
We shall need the following result.
\bl \label{measure-change-1}
 Let  $\varphi: (0,\infty)\rightarrow (0,\infty)$ be a continuous function. For each bounded Borel function $g: S^{n-1}\rightarrow \R$ and convex body $K\in\cK_0^n$,  one has
 \begin{equation}\label{change variable-1}
 \int_{S^{n-1}}g(v)\,d\widetilde{C}_\varphi(K,v)=
 \frac{1}{n}\int_{\partial^\prime K}\big[x\cdot \nu_K(x)\big]\cdot g(\nu_K(x))\cdot \frac{\varphi(|x|)}{|x|^n}\,d\mathcal{H}^{n-1}(x).
 \end{equation}
 \el
\begin{proof}  For each bounded integrable function $f: S^{n-1}\rightarrow \R$,   \cite[(2.31)]{HLYZ} asserts that
\begin{equation}\label{formula 2-31}
 \int_{S^{n-1}}f(u)\rho_K(u)^ndu=\int_{\partial^\prime K}\big[x\cdot \nu_K(x)\big] f(\bar{x}) d\mathcal{H}^{n-1}(x).
\end{equation}  As  $\rho_K$ is positive continuous on $\sphere$ and  $\varphi$ is continuous, one sees that the function $f\cdot \varphi(\rho_K)/\rho_K^n$ is  bounded integrable on $\sphere$. Formula (\ref{formula 2-31}) then implies that  \begin{eqnarray*}
\int_{S^{n-1}}f(u)\varphi(\rho_K(u))\,du
 &=&\int_{\partial^\prime K}\big[x\cdot \nu_K(x)\big]\cdot f(\bar{x}) \cdot \frac{\varphi(\rho_K(\bar{x}))}{\rho_K^n(\bar{x})}\,d\mathcal{H}^{n-1}(x)\\
 &=&\int_{\partial^\prime K}\big[x\cdot \nu_K(x)\big] \cdot f(\bar{x})\cdot \frac{\varphi(|x|)}{|x|^n} \,d\mathcal{H}^{n-1}(x).
 \end{eqnarray*} For each bounded Borel function $g: S^{n-1}\rightarrow \R$, let  $f=g\circ \alpha_K$ which is bounded integrable on $S^{n-1}$.
Then formula (\ref{integral-general curvature measure}) implies  \begin{eqnarray*}
 \int_{S^{n-1}}g(v)\,d\widetilde{C}_\varphi(K,v) &=&
 \frac{1}{n}\int_{S^{n-1}}g(\alpha_K(u))\varphi(\rho_K(u))\,du\\ &=& \frac{1}{n} \int_{\partial^\prime K}\big[x\cdot \nu_K(x)\big] \cdot g(\alpha_K(\bar{x}))\cdot \frac{\varphi(|x|)}{|x|^n} \,d\mathcal{H}^{n-1}(x)\\ &=&\frac{1}{n}  \int_{\partial^\prime K}\big[x\cdot \nu_K(x)\big] \cdot g(\nu_K(x))\cdot \frac{\varphi(|x|)}{|x|^n} \,d\mathcal{H}^{n-1}(x)
 \end{eqnarray*} as desired.  \end{proof}

The dual Orlicz curvature measure has properties similar to those for the $q$-th dual curvature measure. In the following proposition, we will prove some of these properties for the dual Orlicz curvature measure.  \bp\label{measure-check}
Let $K\in\cK_0^n$ and $\varphi: (0,\infty)\rightarrow (0,\infty)$ be a continuous function. The dual Orlicz curvature measure  $\widetilde{C}_\varphi(K,\cdot)$ has the following properties:
\begin{itemize} \item [i)]  $\widetilde{C}_\varphi(K,\cdot)$ is a  Borel measure on $S^{n-1}$;
\item [ii)] $\widetilde{C}_\varphi(K,\cdot)$ is absolutely continuous with respect to the surface area measure $S(K,\cdot)$;
\item [iii)] If the sequence $\{K_i\}_{i=1}^{\infty}\subset \cK_0^n$ converges to $K$ with respect to the Hausdorff metric, then $\widetilde{C}_\varphi(K_i,\cdot)\rightarrow \widetilde{C}_\varphi(K,\cdot)$ weakly.
\end{itemize}
 \ep

 \begin{proof}  i) It is clear that $\widetilde{C}_\varphi(K,\emptyset)=0.$ We only need to prove the countable additivity. Namely, given  a sequence of disjoint Borel sets $\eta_i\subset S^{n-1}$, $i=1,2,\cdots$, with $\eta_i\cap\eta_j=\emptyset$ for $i\neq j$, the following formula holds:
\begin{eqnarray*}
	\widetilde{C}_\varphi(K,\cup_{i=1}^\infty\eta_i)
=\sum_{i=1}^\infty\widetilde{C}_\varphi(K,\eta_i).
\end{eqnarray*}
To this end, it follows from (\ref{formula-general curvature measure}) that for each Borel set $\eta_i\subset \sphere$, one has
$$
 \widetilde{C}_\varphi(K,\eta_i)=\frac{1}{n}\int_{\pmb{\alpha}^*_K(\eta_i)}\varphi(\rho_K(u))du.
$$
By Lemma \ref{properties of maps}, the additivity for Lebesgue integral and the fact that the spherical measure of $\omega_K$ is zero, one has
\begin{eqnarray*}
	\widetilde{C}_\varphi(K,\cup_{i=1}^\infty\eta_i)
	 &=&\frac{1}{n}\int_{\pmb{\alpha}^*_K(\cup_{i=1}^\infty\eta_i)}\varphi(\rho_K(u))du\\
	 &=&\frac{1}{n}\int_{\cup_{i=1}^\infty\pmb{\alpha}^*_K(\eta_i)}\varphi(\rho_K(u))du\\
	&=&\frac{1}{n}\int_{\cup_{i=1}^\infty(\pmb{\alpha}^*_K(\eta_i)\setminus \omega_K)}\varphi(\rho_K(u))du\\
	&=&\frac{1}{n}\sum_{i=1}^\infty\int_{\pmb{\alpha}^*_K(\eta_i)\setminus \omega_K}\varphi(\rho_K(u))du\\
	 &=&\frac{1}{n}\sum_{i=1}^\infty\int_{\pmb{\alpha}^*_K(\eta_i)}\varphi(\rho_K(u))du\\
	&=&\sum_{i=1}^\infty\widetilde{C}_\varphi(K,\eta_i).
\end{eqnarray*}
The countable additivity holds and hence $\widetilde{C}_\varphi(K, \cdot)$ is a Borel measure.

 \vskip 2mm \noindent ii)
 As $K\in\cK_0^n$ and $\varphi$ is continuous, there exists a positive constant $C<\infty$ such that $$[x\cdot\nu_K(x)]\cdot \frac{\varphi(|x|)}{|x|^n}\leq nC$$ for all $x\in\partial K$. Let $\eta \subset S^{n-1}$ be such that $S(K,\eta)=0$ and hence $\mathcal{H}^{n-1}(\nu_K^{-1}(\eta))=0$. It follows from (\ref{change variable-1}) with  $g=\pmb{1}_\eta$  that
\begin{eqnarray*}
 \widetilde{C}_\varphi(K,\eta)
 &=& \frac{1}{n}\int_{\nu_K^{-1}(\eta)}[x\cdot\nu_K(x)]\cdot \frac{\varphi(|x|)}{|x|^n}
 d\mathcal{H}^{n-1}(x)\\
 &\le&C\int_{\nu_K^{-1}(\eta)}d\mathcal{H}^{n-1}(x)\\
 &=&C\cdot \mathcal{H}^{n-1}(\nu_K^{-1}(\eta))\\
 &=&0.
\end{eqnarray*}
That is, $\widetilde{C}_\varphi(K,\cdot)$ is absolutely continuous with respect to $S(K,\cdot)$.

\vskip 2mm \noindent  iii)
 Note that $\rho_{K_i}\rightarrow \rho_{K}$ uniformly and $\alpha_{K_i}\rightarrow \alpha_K$ almost everywhere on $S^{n-1}$ (see Lemma \ref{properties of maps}), due to the convergence $K_i\rightarrow K$ with respect to the Hausdorff metric. Given any continuous function  $g: S^{n-1}\rightarrow \R$,
together with the continuity of $\varphi$, one can find a constant $M>0$ such that for all $i=1,2,\cdots$,
 $$
 |g(\alpha_{K_i})\varphi(\rho_{K_i})|\leq M \ \ \ \ \mathrm{and}\ \ \ |g(\alpha_{K})\varphi(\rho_{K})| \leq M.
 $$
It follows from (\ref{integral-general curvature measure}) and the dominated convergence theorem that
 \begin{eqnarray*}
 	\lim_{i\rightarrow \infty}\int_{S^{n-1}}g(v)\,d\widetilde{C}_\varphi(K_i,v)
 	&=&\lim_{i\rightarrow \infty}\frac{1}{n}\int_{S^{n-1}}g(\alpha_{K_i}(u))\varphi(\rho_{K_i}(u))\,du\\
 	&=&\frac{1}{n}\int_{S^{n-1}}g(\alpha_{K}(u))\varphi(\rho_{K}(u))\,du\\
 	&=&\int_{S^{n-1}}g(v)\,d\widetilde{C}_\varphi(K,v).
 \end{eqnarray*}
 This shows that $\widetilde{C}_\varphi(K_i,\cdot)\rightarrow \widetilde{C}_\varphi(K,\cdot)$ weakly.
 \end{proof}

The following theorem regards to the unique determination of convex bodies by the dual Orlicz curvature measure. When $\varphi(t)=t^q$ with $q<0$, one recovers \cite[Theorem 5.2]{zhao}, which was used to prove the uniqueness of the solution to  the dual Minkowski problem for negative $q$. The techniques used in the proof of \cite[Theorem 5.2]{zhao} (or Theorem \ref{uniqueness-1-1} below) seem not working for strictly increasing function $\varphi$.    \bt\label{uniqueness-1-1}
Suppose that $\varphi: (0, \infty)\rightarrow (0, \infty)$ is a strictly decreasing continuous function. If $K,L\in \cK_0^n$ satisfy that $\widetilde{C}_\varphi(K,\cdot)=\widetilde{C}_\varphi(L,\cdot)$, then $K=L$. \et
\begin{proof} Our proof adopts the beautiful techniques from that of   \cite[Theorem 5.2]{zhao}. First of all, assume that $K, L\in \cK_0^n$ are not dilate to each other. Then there exists a constant $\lambda_0>0$ such that the following sets are nonempty:
\begin{eqnarray*} 
\eta^\prime&=&\{v\in S^{n-1}:h_{K^\prime}(v)>h_L(v)\}, \\ 
\eta &=&\{v\in S^{n-1}:h_{K^\prime}(v)<h_L(v)\},
\\ \eta_{0}&=&\{v\in S^{n-1}:h_{K^\prime}(v)=h_L(v)\},
\end{eqnarray*} where $K^\prime=\lambda_0 K$. It can be easily checked that $\pmb{\alpha}^*_{K}(\omega)=\pmb{\alpha}^*_{K^\prime}(\omega)$ for all $\omega\subset \sphere$.

By \cite[Lemma 5.1 (d)]{zhao}, the set $\pmb{\alpha}^*_{L}(\eta')$ has positive spherical measure. Together with the assumption $\widetilde{C}_\varphi(K,\cdot)=\widetilde{C}_\varphi(L,\cdot)$, one has  $$\widetilde{C}_\varphi(K,\eta^\prime) = \widetilde{C}_\varphi(L,\eta^\prime) 
=\frac{1}{n}\int_{\pmb{\alpha}^*_{L}(\eta^\prime)}\varphi(\rho_L(u))du>0.$$ This further implies that  the set $\pmb{\alpha}^*_{K}(\eta^\prime)$ (and hence $\pmb{\alpha}^*_{K^\prime}(\eta^\prime)$)  has positive spherical measure. Together with \cite[Lemma 5.1 (a) and (c)]{zhao} and the fact that $\varphi$ is strictly decreasing, one gets 
\begin{eqnarray}
\widetilde{C}_\varphi(K,\eta^\prime) =\frac{1}{n}\int_{\pmb{\alpha}^*_{L}(\eta^\prime)}\varphi(\rho_L(u))du\ge\frac{1}{n}\int_{\pmb{\alpha}^*_{K^\prime}(\eta^\prime)}\varphi(\rho_L(u))du >\frac{1}{n}\int_{\pmb{\alpha}^*_{K^\prime}(\eta^\prime)}\varphi(\lambda_0 \rho_{K}(u))du. \label{comparation-1}
 \end{eqnarray} Assume that $\lambda_0\le 1$. Then $\varphi(\rho_K)\le \varphi(\lambda_0 \rho_{K})$ on $\pmb{\alpha}^*_{K}(\eta^\prime)=\pmb{\alpha}^*_{K'}(\eta^\prime)$ as $\varphi$ is strictly decreasing. Thus,  the following inequality holds, which contradicts with (\ref{comparation-1}):  \begin{eqnarray*}
\widetilde{C}_\varphi(K,\eta^\prime)=
\frac{1}{n}\int_{\pmb{\alpha}_K^*(\eta^\prime)}\varphi(\rho_K(u))du\leq \frac{1}{n}\int_{\pmb{\alpha}_{K'}^*(\eta^\prime)}\varphi(\lambda_0 \rho_{K}(u))du.  
\end{eqnarray*} Hence  $\lambda_0>1$ and then  $\varphi(\rho_K)> \varphi(\lambda_0 \rho_{K})$ on $\pmb{\alpha}^*_{K}(\eta)=\pmb{\alpha}^*_{K'}(\eta)$, again due to the fact that $\varphi$ is strictly decreasing. It follows from \cite[Lemma 5.1 (d)]{zhao} that  \begin{eqnarray*}
\widetilde{C}_\varphi(K,\eta)=
\frac{1}{n}\int_{\pmb{\alpha}_K^*(\eta)}\varphi(\rho_K(u))du> \frac{1}{n}\int_{\pmb{\alpha}_{K'}^*(\eta)}\varphi(\lambda_0 \rho_{K}(u))du=\widetilde{C}_\varphi(K^\prime,\eta)>0.  
\end{eqnarray*} Thus the set $\pmb{\alpha}^*_{L}(\eta)$ has positive spherical measure because $$\widetilde{C}_\varphi(K,\eta)=\widetilde{C}_\varphi(L,\eta)=\frac{1}{n}\int_{\pmb{\alpha}^*_{L}(\eta)}\varphi(\rho_L(u))du >0.  $$  Again by \cite[Lemma 5.1]{zhao} and the fact that $\varphi$ is strictly decreasing,  one has 
\begin{eqnarray*} \nonumber
\widetilde{C}_\varphi(K,\eta)   
&=&\frac{1}{n}\int_{\pmb{\alpha}^*_{L}(\eta)}\varphi(\rho_L(u))du\\  
&<&\frac{1}{n}\int_{\pmb{\alpha}^*_{L}(\eta)}\varphi(\lambda_0 \rho_{K}(u))du\\ 
&<&\frac{1}{n}\int_{\pmb{\alpha}^*_{K^\prime}(\eta)}\varphi(\rho_{K}(u))du\\
&=&\frac{1}{n}\int_{\pmb{\alpha}^*_{K}(\eta)}\varphi(\rho_{K}(u))du\\
&=&\widetilde{C}_\varphi(K,\eta), 
 \end{eqnarray*} which is impossible.  This concludes that $K$ and $L$ must be dilates of each other. 
 
 Secondly, without loss of generality, assume that $K\neq L$ but $L=\lambda K$ for some constant $\lambda>1$ (otherwise, switching the roles of $K$ and $L$). It follows from $\widetilde{C}_\varphi(K,\cdot)=\widetilde{C}_\varphi(L,\cdot)$ that 
\begin{eqnarray*} 
\frac{1}{n}\int_{\sphere} \varphi(\rho_K(u))du = \frac{1}{n}\int_{\sphere} \varphi(\rho_L(u))du = \frac{1}{n}\int_{\sphere}\varphi(\lambda \rho_{ K}(u))du< \frac{1}{n}\int_{\sphere }\varphi(\rho_K(u))du,
 \end{eqnarray*} where the last inequality follows from the fact that $\varphi$ is strictly decreasing on $(0,\infty)$. This is a contradiction and hence $K=L$.  
\end{proof}

\section{A variational formula for the dual Orlicz-quermassintegral}\label{variation-2--2--1-1}
Throughout this section, let $\Omega\subset S^{n-1}$ be a closed set that is not contained in any closed hemisphere, and let $h_0, g, \rho_0: \Omega\rightarrow (0,\infty)$ be continuous functions. Denote by $\eta_0 =\eta_{\langle \rho_0\rangle} \subset S^{n-1}$ the spherical measure zero set consisting of the complement of the regular normal vectors of $\langle \rho_0\rangle$.  For each $t\in (-\delta, \delta)$ with $\delta>0$ a fixed constant, let $o(t,\cdot): \Omega\rightarrow (0,\infty)$ be a continuous function such that $\lim_{t\rightarrow 0} o(t,\cdot)/t=0$ uniformly on $\Omega$. Define continuous functions $h_t, \rho_t : \Omega\rightarrow (0,\infty)$  for each $t\in (-\delta, \delta),$ respectively, by  \begin{eqnarray*} \log h_t(v)&=&\log h_0(v)+t g(v)+o(t,v), \\ \log \rho_t(v)&=&\log \rho_0(v)+t g(v)+o(t,v), \ \ \ \mathrm{for\ all}\ v\in \Omega.
\end{eqnarray*} As in Section \ref{section 2}, we use $[h_t]$ and $\langle \rho_t\rangle$ to denote the Wulff shape associated to $h_t$ and the convex hull generated by $\rho_t$, respectively.

Assume that functions $\varphi$ and $\phi$ satisfy the following assumptions:
 \begin{itemize} \item[A1)]  $\phi: (0, \infty)\rightarrow (0, \infty)$ is a strictly decreasing function with $$\lim_{t\rightarrow  0^+}\phi(t)=\infty \ \ \ \mathrm{and} \ \ \ \lim_{t\rightarrow  \infty}\phi(t)=0;$$
  \item [A2)]   $\phi'$, the derivative of $\phi$, exists and is strictly negative on $(0, \infty)$; \item [A3)]  $\varphi(t)=-\phi'(t)t: (0, \infty)\rightarrow (0, \infty)$ is continuous; hence    \begin{equation*}  \phi(t)=\int_t^\infty \frac{\varphi(s)}{s}\,ds. \end{equation*} \end{itemize}   The assumptions  $\lim_{t\rightarrow  0^+}\phi(t)=\infty$  and $\lim_{t\rightarrow  \infty}\phi(t)=0$ are mainly for convenience (especially in Section \ref{solution-dual-Orlicz}). Our results may still work for general strictly decreasing function $\phi$. 

The following lemma is \cite[Lemma 4.1]{HLYZ}, which is essential in the proof of the variational formula for the dual Orlicz-quermassintegral. Let $\Omega\subset S^{n-1}$ be a closed set that is not contained in any closed hemisphere, and  \begin{eqnarray*}  \log \rho_t(v)=\log \rho_0(v)+t g(v)+o(t,v), \ \ \ \mathrm{for\ all}\ v\in \Omega.
\end{eqnarray*}

\bl\label{bounded-lemma-1} For all $v\in S^{n-1}\setminus \eta_0$, one has
 $$ \lim_{t\rightarrow 0}\frac{\log h_{\langle \rho_t\rangle}(v)-\log h_{\langle \rho_0\rangle}(v)}{t}=g(\alpha^*_{\langle \rho_0\rangle}(v)).$$
   Moreover, there exist $\delta>0$ and $M>0$ such that
  $$
 |\log h_{\langle \rho_t\rangle}(v)-\log h_{\langle \rho_0\rangle}(v)|\le M|t|,
 $$
 for all $v\in S^{n-1}$ and all $t\in (-\delta,\delta)$.
 \el

 We now establish the asymptotic behaviour of $\phi(h^{-1}_{\langle \rho_t\rangle})$ as $t\rightarrow 0$ based on  Lemma \ref{bounded-lemma-1}. When $\varphi(t)=t^q$ ($0\neq q\in \R$), it becomes \cite[Lemma 4.2]{HLYZ}.

  \bl\label{bounded-lemma-2}  Suppose that $\varphi$ and $\phi$ satisfy conditions A1)-A3). For all $v\in S^{n-1}\backslash \eta_0$, one has
 \begin{equation*}
 \lim_{t\rightarrow 0}\frac{\phi(h^{-1}_{\langle \rho_t\rangle}(v))-\phi(h^{-1}_{\langle \rho_0\rangle}(v))}{t}=\varphi(h^{-1}_{\langle \rho_0\rangle}(v))g(\alpha^*_{\langle \rho_0\rangle}(v)).
 \end{equation*}
 Moreover, there exist $\delta>0$ and $M>0$ such that
  \begin{equation}\label{bounded-1-1-1}
 |\phi(h^{-1}_{\langle \rho_t\rangle}(v))-\phi(h^{-1}_{\langle \rho_0\rangle}(v))|\le M|t|,
 \end{equation}
 for all $v\in S^{n-1}$ and all $t\in (-\delta,\delta)$.
 \el

 \begin{proof} Recall that $\varphi(t)=-\phi'(t)t$.   Lemma \ref{bounded-lemma-1} and the chain rule for derivative yield that for all $v\in S^{n-1}\backslash \eta_0$,
 \begin{eqnarray*}
 \lim_{t\rightarrow 0}\frac{\phi(h^{-1}_{\langle \rho_t\rangle}(v))-\phi(h^{-1}_{\langle \rho_0\rangle}(v))}{t}&=&\lim_{t\rightarrow 0}\frac{\phi(\exp(-\log h_{\langle \rho_t\rangle}(v)))-\phi(\exp(-\log h_{\langle \rho_0\rangle}(v)))}{t} \\
 &=&-\frac{\phi^\prime(h^{-1}_{\langle \rho_0\rangle}(v))}{h_{\langle \rho_0\rangle}(v)}\cdot \lim_{t\rightarrow 0} \frac{\log h_{\langle \rho_t\rangle}(v)-\log h_{\langle \rho_0\rangle}(v)}{t}\\
 &=&\varphi(h^{-1}_{\langle \rho_0\rangle}(v))\cdot g(\alpha^*_{\langle \rho_0\rangle}(v)).
 \end{eqnarray*} It follows from the uniform convergence of $h_{\langle \rho_t\rangle}\rightarrow h_{\langle \rho_0\rangle}$ that $\{h_{\langle \rho_t\rangle}\}$ is uniformly bounded from both sides, namely,   there exist constants $m_0, m_1, \delta^\prime>0$, such that,  for each $t\in(-\delta^\prime,\delta^\prime)$,  \begin{equation} \label{uniform-bounded-0316} m_0<h_{\langle \rho_t\rangle}<m_1\ \ \ \mathrm{on}\ \ \   S^{n-1}.\end{equation}  Due to the continuity of $\phi$  on $[1/m_1, 1/m_0]$, a constant $M_1>1$ can be found so that
 $$
 0<\phi(h^{-1}_{\langle \rho_t\rangle})/\phi(h^{-1}_{\langle \rho_0\rangle})<M_1 \ \ \text{on}\ \ S^{n-1}.
 $$ Note that  $|s-1|\le M_1 |\log s|$ for $s\in (0,M_1)$ (see e.g. \cite[p.362]{HLYZ}). With $s=\phi(h^{-1}_{\langle \rho_t\rangle})/\phi(h^{-1}_{\langle \rho_0\rangle})$, the following inequality holds on $\sphere$: 
 \begin{equation}\label{bounded-2}
|\phi(h^{-1}_{\langle \rho_t\rangle})-\phi(h^{-1}_{\langle \rho_0\rangle})|\le \phi(h^{-1}_{\langle \rho_0\rangle})\cdot M_1 \cdot |\log \phi(h^{-1}_{\langle \rho_t\rangle})-\log \phi(h^{-1}_{\langle \rho_0\rangle})|.
 \end{equation}

 Note that  on $\sphere$,  $\log h_{\langle \rho_0\rangle}\! \in (\log m_0, \log m_1)$ and $\log h_{\langle \rho_t\rangle}\!\in (\log m_0, \log m_1)$ for each $t\in(-\delta^\prime,\delta^\prime)$. On $(\log (m_0/2),  \log (2 m_1))$, the function $\log \phi(\exp(-s))$ is clearly continuous and differentiable. Due to the continuity of $\varphi$ and $\phi$ on $[1/m_1, 1/m_0]$, there is a constant $M_2>0$ such that for all $s\in (\log m_0, \log m_1)$,   $$\big[\log \phi(\exp(-s))\big]'=\frac{ \varphi(\exp(-s))}{ \phi(\exp(-s))}\leq M_2.$$ It follows from the mean value theorem that for all $s, s'\in (\log m_0, \log m_1)$,   $$ \big|\log \phi(\exp(-s))-\log \phi(\exp(-s'))\big| \leq M_2 \big|s-s'\big|.$$ In particular, for all $v\in \sphere$,  with $s=\log h_{\langle \rho_t\rangle}(v)$ and $s'=\log h_{\langle \rho_0\rangle}(v)$, one has  \begin{eqnarray} \big|\log \phi(h^{-1}_{\langle \rho_t\rangle}(v))-\log \phi(h^{-1}_{\langle \rho_0\rangle}(v))\big|
\leq M_2 \big|\log h_{\langle \rho_t\rangle}(v)-\log h_{\langle \rho_0\rangle} (v)\big|. \label{bounded-3}
 \end{eqnarray} The desired inequality (\ref{bounded-1-1-1}) follows immediately from  (\ref{bounded-2}), (\ref{bounded-3}) and Lemma \ref{bounded-lemma-1}.
\end{proof}

The following theorem provides a variational formula for the dual Orlicz-quermassintegral,  which is the key to solve the dual Orlicz-Minkowski problem based on the method of Lagrange multipliers. When $\varphi(t)=t^q$ with $0\neq q\in\R$, one gets the result in \cite[Lemma 4.5]{HLYZ}. Let $\Omega\subset S^{n-1}$ be a closed set that is not contained in any closed hemisphere, and  \begin{eqnarray*} \log h_t(v) = \log h_0(v)+t g(v)+o(t,v),  \ \ \ \mathrm{for\ all}\ v\in \Omega.
\end{eqnarray*}

 \bt\label{variational-for-decreaing-1}  Suppose that $\varphi$ and $\phi$ satisfy conditions A1)-A3). Given two continuous functions  $h_0: \Omega\rightarrow (0,\infty)$ and $g: \Omega\rightarrow \R$,  one has  \begin{equation}\label{variational formula-2}
 \lim_{t\rightarrow 0}\frac{\widetilde{V}_\phi([h_t])-\widetilde{V}_\phi([h_0])}{t}=- \int_\Omega g(u)\,d\widetilde{C}_\varphi([h_0],u),
 \end{equation} where $[h_t]$ is the family of Wulff shapes associated to $h_t$.   Moreover,
  \begin{equation}\label{variational formula-3}
 \frac{d}{dt}\log\widetilde{V}_\phi([h_t])\bigg |_{t=0}= \frac{-1}{\widetilde{V}_\phi([h_0])}\int_\Omega g(u)\, d\widetilde{C}_\varphi([h_0],u).
 \end{equation} \et
 \begin{proof} Let $\rho_0: \Omega\rightarrow (0,\infty)$ be a continuous function and  $\langle \rho_0\rangle=\text{conv}\{\rho_0(u)u: u\in \Omega\}$. It has been proved that $g$ can be extended to a continuous function $\hat{g}: S^{n-1}\rightarrow \R$  (see \cite[p.364]{HLYZ}),  such that, for all $v\in S^{n-1}\setminus \eta_0$, \begin{equation}\label{extended function}
 g(\alpha_{\langle \rho_0\rangle^*}(v))=(\hat{g}\textbf{1}_\Omega)(\alpha_{\langle \rho_0\rangle^*}(v)).
 \end{equation} Recall that $\log \rho_t(v)=\log \rho_0(v)+t g(v)+o(t,v)$  for all $v\in \Omega.$ It follows from (\ref{formula-volume}), Lemma  \ref{bounded-lemma-2}, (\ref{extended function}), Lemmas  \ref{relation-radial Gauss image} and \ref{measure-change} that
\begin{eqnarray}\nonumber
\lim_{t\rightarrow 0}\frac{\widetilde{V}_\phi(\langle \rho_t\rangle^*)-\widetilde{V}_\phi(\langle \rho_0\rangle^*)}{t}
&=&\lim_{t\rightarrow 0}\frac{1}{n}\int_{S^{n-1}}\frac{\phi(\rho_{\langle \rho_t\rangle^*}(v))-\phi(\rho_{\langle \rho_0\rangle^*}(v))}{t}\, dv\\ \nonumber
&=&\lim_{t\rightarrow 0}\frac{1}{n}\int_{S^{n-1}}\frac{\phi(h^{-1}_{\langle \rho_t\rangle}(v))-\phi(h^{-1}_{\langle \rho_0\rangle}(v))}{t}\,dv\\ \nonumber
&=&\frac{1}{n}\int_{S^{n-1}\setminus \eta_0}\lim_{t\rightarrow 0}\frac{\phi(h^{-1}_{\langle \rho_t\rangle}(v))-\phi(h^{-1}_{\langle \rho_0\rangle}(v))}{t}\, dv\\ \nonumber
&=&\frac{1}{n}\int_{S^{n-1}\setminus \eta_0} \varphi(h^{-1}_{\langle \rho_0\rangle}(v))g(\alpha^*_{\langle \rho_0\rangle}(v))\, dv\\ \nonumber
&=&\frac{1}{n}\int_{S^{n-1}} (\hat{g}\textbf{1}_\Omega)(\alpha_{\langle \rho_0\rangle^*}(v))\varphi(\rho_{\langle \rho_0\rangle^*}(v))\, dv\\ \nonumber
&=&\int_{S^{n-1}} (\hat{g}\textbf{1}_\Omega)(u)\, d\widetilde{C}_\varphi(\langle \rho_0\rangle^*,u)\\
&=&\int_{\Omega} g(u)\, d\widetilde{C}_\varphi(\langle \rho_0\rangle^*,u).\label{variational}
 \end{eqnarray}
Let  $[h_t]$ be the Wulff shape associated to $h_t$ with  $
\log h_t=\log h_0+t g+ o(t,\cdot)$ and  $\kappa_t=1/h_t$. It follows from Lemma \ref{relation} that  $
[h_t]=\langle \kappa_t\rangle^*$ with $\langle \kappa_t\rangle$ the convex hull generated by $\kappa_t$. Note that $$\log \kappa_t=\log \kappa_0-t g- o(t,\cdot).$$  The desired formula (\ref{variational formula-2}) is an immediate consequence of (\ref{variational}) with $\rho_t$ replaced  by $\kappa_t$.

Formula (\ref{variational formula-3}) easily follows from  (\ref{variational formula-2}) and the chain rule for derivative:  \begin{eqnarray*}
\frac{d}{dt}\log\widetilde{V}_\phi([h_t])\bigg |_{t=0}
= \left( \frac{1}{\widetilde{V}_\phi([h_t])}\cdot \frac{d}{dt} \widetilde{V}_\phi([h_t])\right)\bigg |_{t=0}
=\frac{-1}{\widetilde{V}_\phi([h_0])}\int_\Omega g(u)\,d\widetilde{C}_\varphi([h_0],u).
\end{eqnarray*}
\end{proof}

  Along the same lines to the proofs of  Lemma \ref{bounded-lemma-2} and Theorem \ref{variational-for-decreaing-1}, we can prove a variational formula for the dual Orlicz-quermassintegral with functions $\varphi$ and $\phi$ satisfy the following conditions:   \begin{itemize} \item[B1):]  $\phi: (0, \infty)\rightarrow (0, \infty)$ is a strictly increasing and continuous function with $$\lim_{t\rightarrow  0^+}\phi(t)=0 \ \ \ \mathrm{and} \ \ \ \lim_{t\rightarrow  \infty}\phi(t)=\infty;$$
  \item [B2):]   $\phi'$, the derivative of $\phi$, exists and is strictly positive on $(0, \infty)$; \item [B3):]  $\varphi(t)=\phi'(t)t: (0, \infty)\rightarrow (0, \infty)$ is a continuous function, and hence    \begin{equation*}  \phi(t)=\int_0^t \frac{\varphi(s)}{s}\,ds. \end{equation*}  \end{itemize}    The assumptions  $\lim_{t\rightarrow  0^+}\phi(t)=0$  and $\lim_{t\rightarrow  \infty}\phi(t)=\infty$ are mainly for convenience; our results may still work for general strictly increasing function $\phi$. 

\bt\label{variational-for-increaing-1}  Suppose that $\varphi$ and $\phi$ satisfy conditions B1)-B3). Given two continuous functions  $h_0: \Omega\rightarrow (0,\infty)$ and $g: \Omega\rightarrow \R$,  one has  
$$
 \lim_{t\rightarrow 0}\frac{\widetilde{V}_\phi([h_t])-\widetilde{V}_\phi([h_0])}{t}= \int_\Omega g(u)\,d\widetilde{C}_\varphi([h_0],u),
 $$ where $[h_t]$ is the family of Wulff shapes associated to $h_t$.   Moreover,
  $$
 \frac{d}{dt}\log\widetilde{V}_\phi([h_t])\bigg |_{t=0}= \frac{1}{\widetilde{V}_\phi([h_0])}\int_\Omega g(u)\, d\widetilde{C}_\varphi([h_0],u).
 $$ \et

\section{A solution to the dual Orlicz-Minkowski problem}\label{solution-dual-Orlicz}

In this section, we provide a solution to the following dual Orlicz-Minkowski problem: {\it under what conditions on $\varphi$ and a given nonzero finite Borel measure $\mu$ on $\sphere$,  there exist a constant $\tau>0$ and a convex body $K$ (ideally with the origin in its interior) such that $\mu=\tau \widetilde{C}_\varphi(K,\cdot)$?}  When $\varphi(t)=t^{q}$ with $0\neq q\leq n$, this problem has been investigated in \cite{HLYZ, zhao, zhao-jdg}.

 For a nonzero finite Borel measure $\mu$ on $\sphere$, let  $$|\mu|=\int_{\sphere}\,d\mu.$$ Recall that for each $K\in \cK_0^n$,   $$\widetilde{V}_\varphi(K)=
 \widetilde{C}_\varphi(K,\sphere)
 =\frac{1}{n}\int_{\sphere}\varphi(\rho_K(u))du.
 $$ Our main result is stated in the following theorem.
 \bt \label{solution-dual-Orlicz-main} Suppose that $\varphi$ and $\phi$ satisfy conditions A1)-A3). Let $\mu$ be a nonzero finite Borel measure on $S^{n-1}$. Then the following statements are equivalent:
\begin{itemize} \item [i)] $\mu$ is not concentrated in any closed hemisphere, i.e.,   $$\int _{\sphere} \langle  \xi, \theta\rangle_+ \,d\mu(\theta)>0\ \ \ \mathrm{for\ \  all\ \ }  \xi \in \sphere;$$
\item [ii)]  there exists a convex body $K\in\cK_0^n$, such that, $$\frac{\mu}{|\mu|}=\frac{\widetilde{C}_\varphi(K,\cdot)}{\widetilde{V}_\varphi(K)}. $$ \end{itemize}
 \et
\begin{proof} First, let us prove the easier direction  ii)$\Rightarrow$i). Suppose that $\mu=\tau \widetilde{C}_\varphi(K,\cdot)$ for some convex body $K\in\cK_0^n$ and $$\tau=\frac{|\mu|}{\widetilde{V}_\varphi(K)}>0.$$ Note that $\varphi$ is continuous on $[r_K, R_K]$, where $$0<r_K=\min_{u\in S^{n-1}} \rho_K(u)  \ \ \mathrm{and} \ \ \ R_K=\max_{u\in S^{n-1}} \rho_K(u)<\infty.$$ Hence, there exists a constant $C>0$, such that, for all $x\in \partial K$,  $$\big[x\cdot \nu_K(x)\big]  \cdot \frac{\varphi(|x|)}{|x|^n}\geq n C.$$  The surface area measure $S(K, \cdot)$ of a convex body $K\in\cK^n_0$ is not concentrated in any closed hemisphere, that is, for any $\xi\in S^{n-1}$,
 \begin{equation}\label{surface-measure-Minkowski solution}
 \int_{S^{n-1}}(\xi \cdot u)_+dS(K, u)> 0.
 \end{equation}

Together with Lemma \ref{measure-change-1}, (\ref{surface area measure}),  (\ref{surface-measure-Minkowski solution}) and $\mu=\tau \widetilde{C}_\varphi(K,\cdot)$, one has, for any $\xi\in S^{n-1}$,
\begin{eqnarray*}
\int_{S^{n-1}}(\xi \cdot v)_+\,d\mu(v)
&=&\tau\int_{S^{n-1}}(\xi \cdot v)_+\,d\widetilde{C}_\varphi(K, v)\\
&=&\frac{\tau}{n}\int_{\partial^\prime K}  \big[(\xi\cdot \nu_K(x))_+\big]\cdot \big[x\cdot \nu_K(x)\big] \cdot \frac{\varphi(|x|)}{|x|^n}\,d\mathcal{H}^{n-1}(x)\\
&\ge&\tau C \int_{\partial^\prime K}(\xi\cdot \nu_K(x))_+\,d\mathcal{H}^{n-1}(x)\\
&=&\tau C \int_{S^{n-1}}(\xi \cdot u)_+\,dS(K, u)\\
&>& 0. \end{eqnarray*}
  This implies that $\mu$ is not concentrated in any closed hemisphere.

 \vskip 2mm Now let us prove the direction i)$\Rightarrow$ii). The proof needs several steps. Let $\mu$ be a nonzero finite Borel measure on $S^{n-1}$ such that  $\mu$ is not concentrated in any closed hemisphere.

 \vskip 2mm \noindent {\it Step 1:} if $\{Q_i\}_{i=1}^\infty\subset \cK_0^n$ and $c>0$ is a constant, such that,  $\widetilde{V}_\phi(Q_i)=c$,
then there exists a constant $R>0$ such that $
Q_i^* \subset R \ball. $

 To this end, we assume that there are no such constants $R$ such that $Q_i^*\subset R\ball$. Without loss of generality, assume that $\{Q_i^*\}_{i=1}^\infty$  satisfies $R_{Q_i^*}\rightarrow \infty$ as $i\rightarrow \infty$, where $R_{Q_i^*}$ is the maximum radius of $Q_i^*$, i.e.,  $$
R_{Q_i^*}=\rho_{Q_i^*}(v_i)=\max\big\{\rho_{Q_i^*}(v), v\in S^{n-1}\big\}.
$$ We can further assume that, due to the compactness of $S^{n-1}$,  $\{v_i\}_{i=1}^\infty\subset S^{n-1}$ is a convergent sequence with limit $v_0\in S^{n-1}$, namely, $\lim_{i\rightarrow \infty} v_i=v_0.$ It is obvious that, for all $u\in \sphere$,  \begin{equation}
h_{Q_i^*}(u)\ge (u\cdot v_i)_+\ R_{Q_i^*}.\label{segment--1}
\end{equation}

Note that the spherical measure is not concentrated in any closed hemisphere.  This fact yields that  there exists a constant $c_0>0$, such that,   \begin{equation*}
 \int_{S^{n-1}}(u\cdot v_0)_+\ du\ge c_0. 
 \end{equation*}   For  all integers $j\geq 1$, let $$\Sigma_j(v_0)=\bigg\{u\in S^{n-1}:(u \cdot v_0)_+ >\frac{1}{j}\bigg\}.$$  One can check that  $\Sigma_j(v_0)$ forms an increasing nest of sets:  $\Sigma_j(v_0)\subset \Sigma_{j+1}(v_0)$ for all $j\geq 1$. Moreover,  $\cup_{j=1}^{\infty}\Sigma_j(v_0)=\{u\in S^{n-1}: (u \cdot v_0)_+ >0\}.$ The monotone convergence theorem yields
\begin{eqnarray*} \lim_{j\rightarrow \infty}\int_{\Sigma_j(v_0)}  (u \cdot v_0)_+ \,du =\int_{\cup_{j=1}^{\infty}\Sigma_j(v_0)} (u \cdot v_0)_+ \,du = \int_{S^{n-1}} (u \cdot v_0)_+ \,du \geq c_0.
\end{eqnarray*} Hence,   there exists an integer  $j_0\geq 1$ such that \begin{equation}
\label{sphere integration} \int_{\Sigma_{j_0}(v_0)}\,du\geq \int_{\Sigma_{j_0}(v_0)}  (u \cdot v_0)_+\,
du\geq \frac{c_0}{2}. \end{equation}

For convenience,  let $G: (0,\infty)\rightarrow (0, \infty)$ be the function given by $
G(t)=\phi\left(t^{-1}\right).$ As  $\phi$ is a strictly decreasing function with $\lim_{t\rightarrow \infty}\phi(t)=0$ and $\lim_{t\rightarrow 0^+}\phi(t)=\infty$,  $G$ is a strictly increasing function with $\lim_{t\rightarrow 0^+}G(t)=0$ and $\lim_{t\rightarrow \infty}G(t)=\infty$. Let $M>0$ be a fixed number.  Fatou's Lemma,  (\ref{segment--1}), (\ref{sphere integration}),   and the fact that $G$ is increasing and continuous  imply
\begin{eqnarray*}\nonumber
 \lim_{i\rightarrow \infty}\widetilde{V}_\phi(Q_i) &=&\lim_{i\rightarrow \infty}\frac{1}{n}\int_{S^{n-1}}\phi(\rho_{Q_i}(u))\ du\\  \nonumber
&=&\lim_{i\rightarrow \infty}\frac{1}{n}\int_{S^{n-1}}G(h_{Q_i^*}(u))\ du\\  \nonumber
&\ge&\liminf_{i\rightarrow \infty}\frac{1}{n}\int_{S^{n-1}}G(R_{Q_i^*}\cdot(u\cdot v_i)_+)\ du\\ \nonumber
&\ge&\liminf_{i\rightarrow \infty}\frac{1}{n}\int_{S^{n-1}}G(M\cdot (u\cdot v_i)_+)\ du\\ \nonumber
&\ge&\frac{1}{n}\int_{S^{n-1}}G(M\cdot (u\cdot v_0)_+)\ du\\ \nonumber
&\ge&\frac{1}{n}\int_{\sum_{j_0}(v_0)}G\left(M/j_0\right)\ du\\ \nonumber
&\ge& G\left(M/j_0\right)\cdot \frac{c_0}{2n}.
 \end{eqnarray*}
Note that $G\left(M/j_0\right)\rightarrow \infty$ as $M\rightarrow\infty$ and  $\widetilde{V}_\phi(Q_i)=c$ for all integers $i$. Together with the above inequalities, a contradiction, i.e., $c\ge \infty$, is obtained. Hence the sequence $\{Q_i^*\}_{i=1}^{\infty}$  is uniformly bounded, namely,   there exists a constant $R>0$ such that $Q_i^*
\subset R \ball.$

\vskip 2mm \noindent {\it Step 2:}  there exists a convex body $Q_0\in\cK_0^n$ such that $\widetilde{V}_\phi(Q_0)=|\mu|$ and
\begin{equation*} 
\Phi(Q_0)=\sup\big\{\Phi(K): \widetilde{V}_\phi(K)=|\mu| \ \text{and}\ K\in \cK_0^n \big\},
\end{equation*}  where  $\Phi: \cK_0^n\rightarrow \R$ is defined  by
 \begin{equation}\label{function-optimal-1}
 \Phi(K)=-\frac{1}{|\mu|}\int_{S^{n-1}}\log h_K(v)d\mu(v).
 \end{equation}

 The proof of this step is almost identical to that of \cite[Lemma 4.2]{zhao}. For completeness, we include a brief proof here.  Let $\{Q_i\}_{i=1}^n\subset \cK_0^n$ be a maximizing sequence such that $\widetilde{V}_\phi(Q_i)=|\mu|$ and
$$
\lim_{i\rightarrow \infty} \Phi(Q_i)=\sup\Big\{\Phi(K): \widetilde{V}_\phi(K)=|\mu| \ \text{and}\ K\in \cK_0^n\Big\}.
$$ It follows from Step 1 that $\{Q_i^*\}_{i=1}^{\infty}$  is uniformly bounded. The Blaschke selection theorem implies the existence of a compact convex set $Q\subset \R^n$ and  a subsequence of $\{Q^*_i\}_{i=1}^{\infty}$, which will not be relabeled, such that,  $
Q_i^*\rightarrow Q.$  Now we show that  $o\in \text{int}(Q)$.  Assume  $o\in \partial Q$ and then there exists $u_0\in S^{n-1}$ such that $$
 \lim_{i\rightarrow \infty} h_{Q_i^*}(u_0)=h_{Q}(u_0)=0.$$ It can be proved that   $
 \mu(\Sigma_{\delta_0}(u_0))>0 $ and  $\rho_{Q_i^*}\rightarrow 0$ uniformly on $\Sigma_{\delta_0} (u_0)$ for some $\delta_0>0$ (see details in  \cite[Lemma 4.2]{zhao}),  where $\Sigma_{\delta_0} (u_0)\subset \sphere$ is given by $$
\Sigma_{\delta_0} (u_0)=\{v\in S^{n-1}: v\cdot u_0>\delta_0\}.$$  Together with  (\ref{function-optimal-1}),  $\widetilde{V}_\phi(Q_i)=|\mu|$ and $Q_i^* \subset R \ball$ for all $i$,  one can get
\begin{eqnarray*}
\Phi(Q_i)   \leq \frac{1}{|\mu|}\int_{\Sigma_{\delta_0}(u_0)}\log \rho_{Q_i^*}(v)d\mu(v)
+\log R
 \end{eqnarray*} and hence  $\lim_{i\rightarrow \infty} \Phi(Q_i) =-\infty$, a contradiction.  In conclusion,   $o\in \text{int}(Q)$ and $Q\in \cK_0^n$, which implies  $Q_0=Q^*\in \cK_0^n$.  Moreover $Q_i\rightarrow Q_0$. The desired claim in Step 2 follows immediately by the continuity of  $\Phi(\cdot)$ and $\widetilde{V}_\phi(\cdot)$ (see Lemma \ref{continuity-dual-qu-1}) on $\cK_0^n$. That is, $
\widetilde{V}_\phi(Q_0)=\lim_{i\rightarrow \infty}\widetilde{V}_\phi(Q_i)=|\mu| $ and
$$\Phi(Q_0)=\lim_{i\rightarrow \infty}\Phi(Q_i)=\sup\Big\{\Phi(K): \widetilde{V}_\phi(K)=|\mu| \ \text{and}\ K\in \cK_0^n\Big\}.$$

\noindent {\it Step 3:} the convex body $Q_0$ found in Step 2  is a solution of the dual Orlicz-Minkowski problem, that is,
  $$\frac{\mu}{|\mu|}=\frac{\widetilde{C}_\varphi(Q_0,\cdot)}{\widetilde{V}_\varphi(Q_0)}. $$

 To this end, consider the following optimization problem:
 \begin{equation}\label{max-problem-2}
\sup\Big\{\Phi(f): \widetilde{V}_\phi([f])=|\mu| \ \text{for}\ f\in C^+(S^{n-1})\Big\},
 \end{equation}  where  the functional $\Phi: C^+(S^{n-1})\rightarrow \R$ given by
 \begin{equation}\label{max-function}
 \Phi(f)=-\frac{1}{|\mu|}\int_{S^{n-1}}\log f(v)d\mu(v).
 \end{equation}  
 From the definition of the Wulff shape,  for all $f\in C^+(S^{n-1})$, one has
$$
\Phi(f)\le \Phi(h_{[f]})\  \ \text{and} \ \  \widetilde{V}_\phi([f])=\widetilde{V}_\phi(h_{[f]}).
$$ Hence, it is enough to find maximizers for the optimization problem (\ref{max-problem-2}) among support functions of convex bodies in $\cK_0^n$. Step 2 implies that $h_{Q_0}$ is a maximizer to (\ref{max-problem-2}).

Let $g\in C(S^{n-1})$ be an arbitrary but fixed continuous function and let $\delta>0$ be a small enough constant.  Let $
h_t(v)=h_{Q_0}(v) e^{tg(v)}$ for $t\in (-\delta, \delta)$ and $v\in S^{n-1}.$ Define the following functional on $C^+(\sphere)$  $$\mathcal{L}(t, \tau)=\Phi(h_t)-\tau (\log \widetilde{V}_\phi([h_t]) -\log |\mu|).$$ As $h_{Q_0}$ is a maximizer to (\ref{max-problem-2}), it follows from the method of Lagrange multipliers that  $h_{Q_0}$ must satisfy the following equation: 
 $$
\frac{\partial}{\partial t}\mathcal{L}(t, \tau)\bigg|_{t=0}=0.
$$
By (\ref{variational formula-3}), (\ref{max-function}) and $\widetilde{V}_\phi(Q_0)=|\mu|$, one gets
\begin{eqnarray*}
0&=&\frac{\partial}{\partial t}\left(\left[-\frac{1}{|\mu|}\int_{S^{n-1}}[\log h_{Q_0}(v)+t g(v)]\,d\mu(v)\right]-\tau \log \widetilde{V}_\phi([h_t]) +\tau\log |\mu|\right)\bigg|_{t=0}\\
&=&-\frac{1}{|\mu|}\int_{S^{n-1}}g(v)\,d\mu(v)
+\frac{\tau}{\widetilde{V}_\phi(Q_0)}\int_{S^{n-1}}g(v)\,d\widetilde{C}_\varphi(Q_0, v)\\
&=&\frac{1}{|\mu|}\left(-\int_{S^{n-1}}g(v)\,d\mu(v)
+ \tau \int_{S^{n-1}}g(v)\,d\widetilde{C}_\varphi(Q_0, v)\right).
 \end{eqnarray*} That is,   for all $g\in C(\sphere)$,
\begin{equation}\label{equation-for-arbitary-g}
\int_{S^{n-1}}g(v)\,d\mu(v)
= \tau \int_{S^{n-1}}g(v)\,d\widetilde{C}_\varphi(Q_0, v),
\end{equation} and hence  $\mu=\tau \widetilde{C}_\varphi(Q_0,\cdot).$ The constant $\tau$ can be easily calculated by (\ref{equation-for-arbitary-g}) with $g=1$: \begin{equation*}
|\mu|=\int_{S^{n-1}}\,d\mu(v)
= \tau \int_{S^{n-1}} \,d\widetilde{C}_\varphi(Q_0, v)=\tau \widetilde{V}_\varphi(Q_0),
\end{equation*} which implies $$\tau=\frac{|\mu|}{\widetilde{V}_\varphi(Q_0)}.$$ Putting the constant $\tau$ into  $\mu=\tau \widetilde{C}_\varphi(Q_0,\cdot),$ one gets $$\frac{\mu}{|\mu|}=\frac{\widetilde{C}_\varphi(Q_0,\cdot)}{\widetilde{V}_\varphi(Q_0)} $$ and hence $Q_0$ is a solution to the dual Orlicz-Minkowski problem.      \end{proof}

  Let $\varphi(t)=t^q$ for $q<0$ and hence $\phi(t)=-t^q/q$, which satisfy conditions A1)-A3). By Theorem \ref{solution-dual-Orlicz-main},  one can get the following solution to the dual Minkowski problem for negative $q$, which has been recently proved in \cite{zhao} by Zhao.
 \bc
Suppose that $q<0$ and $\mu$ is a finite nonzero Borel measure on $\sphere$. The following statements are equivalent: \begin{itemize} \item [i)] $\mu$ is not concentrated in any closed hemisphere, i.e.,   $$\int _{\sphere} \langle \xi, \theta\rangle_+ \,d\mu(\theta)>0\ \ \ \mathrm{for\ \  all\ \ } \xi\in \sphere;$$
\item [ii)] there exists  a convex body $K\in\cK_0^n$, such that, $\mu=\widetilde{C}_q(K, \cdot)$. \end{itemize}
 \ec

Note that Zhao  \cite[Theorem 5.2]{zhao}  also proved that the solution to  the dual Minkowski problem for negative $q$ must be unique. Due to lack of homogeneity of the function $\varphi$, to prove the uniqueness of the solutions in the Orlicz setting seems very intractable.  We would like to mention that there are no arguments regarding the  uniqueness of solutions to the Orlicz-Minkowski problem \cite{HS-JFA-1, huanghe2012, liaijun2014}.

\vskip 2mm \noindent {\bf Acknowledgments.} The first author is
supported by AARMS, NSERC, NSFC (No.\ 11501185) and the Doctor
Starting Foundation of Hubei University for Nationalities (No.\
MY2014B001). The third author is supported by a NSERC
grant.

\vskip 2mm \noindent Baocheng Zhu, \ \ \ {\small \tt zhubaocheng814@163.com}\\
{ \em 1. Department of Mathematics,
 Hubei University for Nationalities,
 Enshi, Hubei, China 445000}\\
{  \em 2.\  Department of Mathematics and Statistics,   Memorial University of Newfoundland,
   St.\ John's, Newfoundland, Canada A1C 5S7 }

\vskip 2mm \noindent Sudan Xing, \ \ \ {\small \tt sudanxing@gmail.com}\\
{ \em Department of Mathematics and Statistics,   Memorial University of Newfoundland,
   St.\ John's, Newfoundland, Canada A1C 5S7 }

\vskip 2mm \noindent Deping Ye, \ \ \ {\small \tt deping.ye@mun.ca}\\
{ \em Department of Mathematics and Statistics,
   Memorial University of Newfoundland,
   St.\ John's, Newfoundland, Canada A1C 5S7 }

\end{document}